# Waveform Transmission Method,

# A New Waveform-relaxation Based Algorithm

# to Solve Ordinary Differential Equations in Parallel

Fei Wei

Huazhong Yang

Department of Electronic Engineering, Tsinghua University, Beijing, China

## Abstract

Waveform Relaxation method (WR) is a distributed algorithm to solve Ordinary Differential Equations (ODEs). In this paper, we propose a new distributed algorithm, named Waveform Transmission Method (WTM), by virtually inserting waveform transmission lines into the dynamical system to achieve distributed computing of ODEs. WTM is convergent to solve linear SPD ODEs.

1. Introduction

   Waveform-relaxation method was a daring attempt to solve ODEs extracted from circuit by the transient analysis [1]. WR makes use of a piece of waveform to iterate among the circuits, instead of a single value at one time point.
   Virtual Transmission Method (VTM) is a new distributed algorithm to solve sparse linear systems [3]. Its physical background is microwave network and lossless transmission line. If we replace each unknown of VTM by a piece of waveform, VTM is expanded into a WR based algorithm, called Waveform Transmission Method (WTM).
   By marrying WR with VTM, WTM is born. WTM is a distributed algorithm to solve sparse ODEs. The waveforms of WTM are transferred in a bidirected way, similar to Gauss-Jacobi..
   This paper is organized as follows. In Section 2 we define the Waveform Transmission Line (WTL). In Section 3 we describe the basic steps of WTM. Section 4 gives a simple example. Finally we conclude this paper in Section 5.

## 2. Waveform Transmission Line (WTL)

Waveform Transmission Line (WTL) is extended from the virtual transmission line [4]. The mathematical description of the virtual transmission line is shown in (2.1). The diagram of the virtual transmission line is shown in Fig. 1. (2.1) is called Transmission Delay Equations.

$$\begin{cases} U_1(p) + ZI_1(p) = U_2(p-\rho) - ZI_2(p-\rho) \\ U_2(p) + ZI_2(p) = U_1(p-\rho) - ZI_1(p-\rho) \end{cases} \quad (2.1)$$

where $U_1(p)$ and $U_2(p)$ represent the potential of the two ports of the virtual transmission line, respectively, while $I_1(p)$ and $I_2(p)$ represent the outflow current. $p$ is the virtual time variable and $\rho$ is the virtual transmission delay. $Z$ is the characteristic impedance, which is positive.

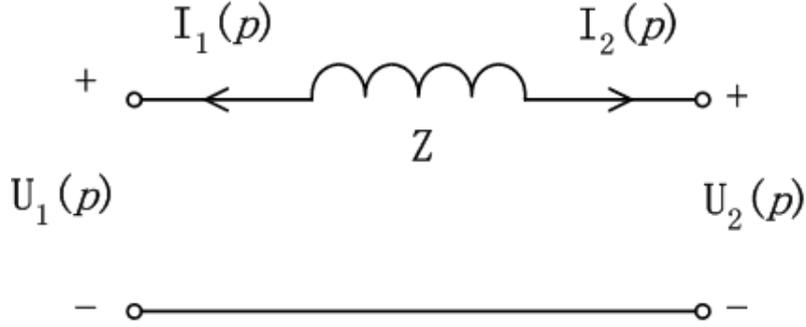

Figure 1. The diagram of the virtual transmission line.

If we replace each unknown of the transmission line by a piece of waveform, we get the Waveform Transmission Line (WTL). The mathematical description of the WTL is shown in (2.2). The diagram of the WTL is shown in Fig. 2. (2.2) is called Waveform Transmission Delay Equations.

$$\begin{cases} U_1(t,p) + Z(t)I_1(t,p) = U_2(t,p-\rho) - Z(t)I_2(t,p-\rho) \\ U_2(t,p) + Z(t)I_2(t,p) = U_1(t,p-\rho) - Z(t)I_1(t,p-\rho) \end{cases} \quad (2.2)$$

$$t \in [T_1, T_2]$$

here $t$ is the physical time variable. $p$ is the virtual time variable. $\rho$ is the virtual transmission delay. $U_1(t,p)$ and $U_2(t,p)$, $t \in [T_1, T_2]$, represent the potential waveforms of the two ports of the waveform transmission line at the virtual time $p$, respectively. $I_1(t,p)$ and $I_2(t,p)$, $t \in [T_1, T_2]$, represent the outflow current waveforms. $Z(t)$, $t \in [T_1, T_2]$, is the characteristic impedance waveform, which should be positive, i.e. $Z(t) > 0$, $t \in [T_1, T_2]$. $Z(t)$ could be considered as a preconditioner for WTM.

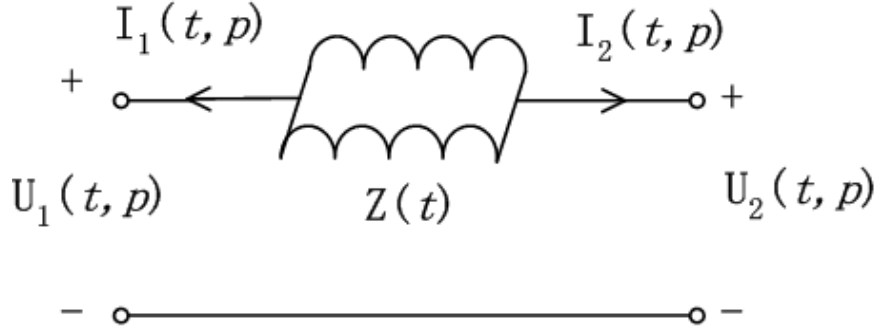

Figure 2. The diagram of the Waveform Transmission Line

If the delay $\rho$ of all the WTLs are same, then (2.2) could be simplified into a discrete iterative form, as shown in (2.3). Consequently, (2.2) is the continuous iterative form of (2.3).

$$\begin{cases} U_1^k(t) + Z(t)I_1^k(t) = U_2^{k-1}(t) - Z(t)I_2^{k-1}(t) \\ U_2^k(t) + Z(t)I_2^k(t) = U_1^{k-1}(t) - Z(t)I_1^{k-1}(t) \end{cases} \quad (2.3)$$
$$t \in [T_1, T_2]$$

3. Waveform Transmission Method

The mathematic description of the linear ODEs is (3.1). If both **C** and **A** are Symmetric Positive Definite (SPD), this kind of ODEs is called SPD ODEs. In this paper, we mainly focus on how to solve SPD ODEs.

$$\mathbf{C}\frac{d\mathbf{x}(t)}{dt} + \mathbf{A} \cdot \mathbf{x}(t) = \mathbf{b}, \quad \mathbf{x}(0) = \mathbf{x}_0 \quad (3.1)$$

The basic step of WTM is similar to VTM [3, 4].
1. Map the ODEs into an electric graph, which contains a weighted graph of **C**, a weighted graph of **A**, and **b** as the vertex source vector.
2. Set the vertex splitting boundary, and perform the Electric Vertex Splitting (EVS) [3, 4].
    2.1 Each vertex on boundary is split into a pair of twin vertices.
    2.2 The weighted graphs of **C** and **A** are electrically split along the boundary, as well as the vertex source **b.**
    2.3 Add inflow currents into the twin vertices. As the result, the original electric graph for ODEs is split into *n* subgraghs.
3. Add one WTL between each pair of twin vertices.
4. Set the time window [$T_1$, $T_2$] for the waveform. Set the initial waveforms for

the potentials and inflow currents of the twin vertices.
5. Locate each subgragh on a processor, and perform the distributed iteration on $n$ processors.

Generally, we point that, if the electric graph of an SPD ODEs (3.1) is partitioned into $n$ subgraghs, and all these subgraghs are SNND, then for positive characteristic impedance waveforms of the waveform transmission lines, WTM converges to the solution of the original system. This conclusion is similar to the convergence theory of VTM [3].

4. Example

Fig 3A is a simple capacity-resistor network, whose mathematic description is shown in (4.1).

$$\begin{cases} C\dfrac{du(t)}{dt} + G \cdot u(t) = b \\ u(0) = u_0 \end{cases} \quad (4.1)$$

here $C = 3$, $G = 1.5$, $u_0 = 0$.

Then, we split this network into two subgraghs by EVS, and insert a zero resistor, as shown in Fig 3B. (4.2) and (4.3) are the mathematic descriptions for these two subgraghs, respectively.

$$\begin{cases} C_1\dfrac{du_1(t)}{dt} + G_1 \cdot u_1(t) = b_1 + i_1(t) \\ u_1(0) = u_0 \end{cases} \quad (4.2)$$

here $C_1 = 1$, $G_1 = 0.5$.

$$\begin{cases} C_2\dfrac{du_2(t)}{dt} + G_2 \cdot u_2(t) = b_2 + i_2(t) \\ u_2(0) = u_0 \end{cases} \quad (4.3)$$

here $C_2 = 2$, $G_2 = 1$.

After that, we replace the zero resistor by one WTL. For simplicity, here we set the characteristic impedance waveform to be a constant waveform, $Z(t) = 1.5$, $t \in [T_1, T_2]$, $T_1 = 0$, $T_2 = 1.0$.

$$\begin{cases} u_1^k(t) + Z(t) \cdot i_1^k(t) = u_2^{k-1}(t) - Z(t) \cdot i_2^{k-1}(t) \\ u_2^k(t) + Z(t) \cdot i_2^k(t) = u_1^{k-1}(t) - Z(t) \cdot i_1^{k-1}(t) \end{cases} \quad (4.4)$$
$$t \in [T_1, T_2]$$

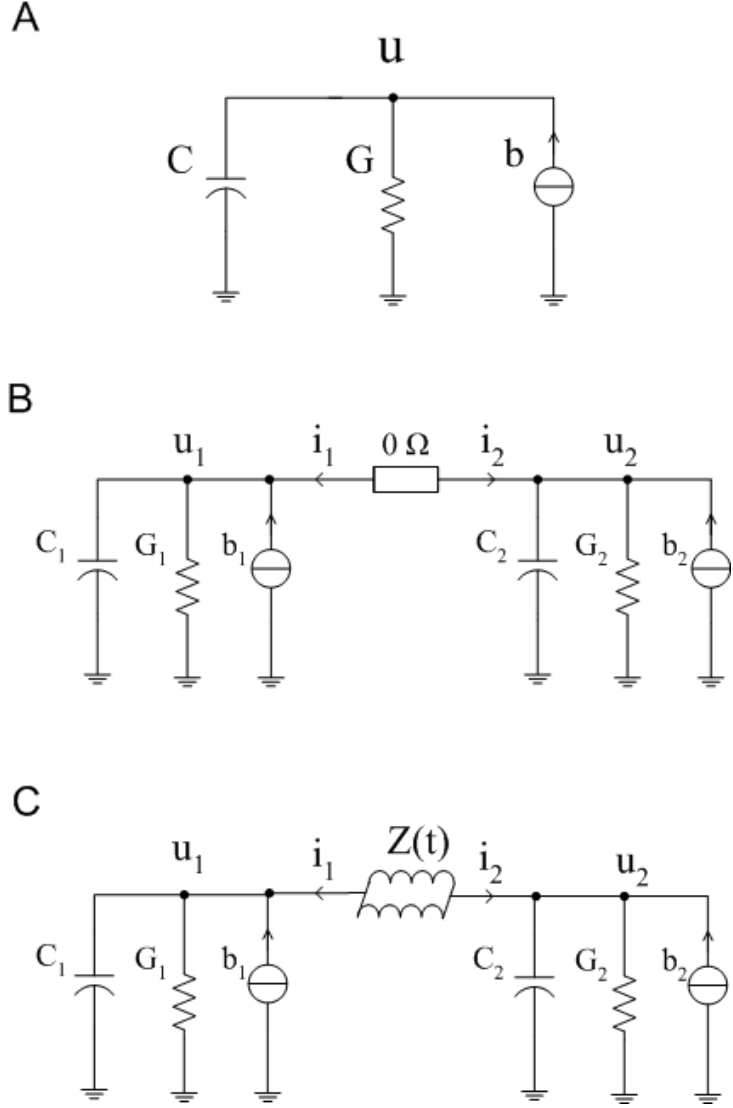

Figure 3.   Illustration of the Electric Vertex Splitting. (A) The original capacity-resistor network. (B) Electrically split the network, and a zero resistor is inserted between them. C = $C_1+C_2$. G = $G_1+G_2$. b = $b_1+b_2$. (C) Replace the zero resistor by a waveform transmission line.

Later, combine (4.2) and (4.4), we get the description for subgragh 1, as below:

$$\begin{cases} C_1 \dfrac{du_1^k(t)}{dt} + G_1 \cdot u_1^k(t) = b_1 + i_1^k(t) \\ u_1^k(t) + Z(t) \cdot i_1^k(t) = u_2^{k-1}(t) - Z(t) \cdot i_2^{k-1}(t) \\ \qquad t \in [T_1, T_2] \end{cases} \quad (4.5)$$

(4.5) is then simplified into (4.6):

$$\begin{cases} C_1 \dfrac{du_1^k(t)}{dt} + \left(G_1 + Z^{-1}(t)\right) \cdot u_1^k = Z^{-1}(t) \cdot u_2^{k-1}(t) - i_2^{k-1}(t) + b_1 \\ i_1^k(t) = -Z^{-1}(t) \cdot u_1^k(t) + Z^{-1}(t) \cdot u_2^{k-1}(t) - i_2^{k-1}(t) \\ t \in [T_1, T_2] \end{cases} \quad (4.6)$$

Similarly, we combine (4.3) and (4.4), and get the description for subgragh 2, as below:

$$\begin{cases} C_2 \dfrac{du_2^k(t)}{dt} + G_2 \cdot u_2^k(t) = b_2 + i_2^k(t) \\ u_2^k(t) + Z(t) \cdot i_2^k(t) = u_1^{k-1}(t) - Z(t) \cdot i_1^{k-1}(t) \\ t \in [T_1, T_2] \end{cases} \quad (4.7)$$

Reformat (4.7) into (4.8):

$$\begin{cases} C_2 \dfrac{du_2^k(t)}{dt} + \left(G_2 + Z^{-1}(t)\right) \cdot u_2^k = Z^{-1}(t) \cdot u_1^{k-1}(t) - i_1^{k-1}(t) + b_2 \\ i_2^k(t) = -Z^{-1}(t) \cdot u_2^k(t) + Z^{-1}(t) \cdot u_1^{k-1}(t) - i_1^{k-1}(t) \\ t \in [T_1, T_2] \end{cases} \quad (4.8)$$

Finally, we locate (4.6) on processor 1, and (4.8) on processor 2, then do the distributed iteration. Numerical experiments show that WTM is convergent. To illustrate the error of WTM after $k$ iterations, we define the max error of waveform, as in (4.9) and (4.10). Fig 4 gives the convergent curve for this example.

$$max\_err\_u_1^k = Max\left(u_1^k(t) - u(t)\right), t \in [T_1, T_2] \quad (4.9)$$

$$max\_err\_u_2^k = Max\left(u_2^k(t) - u(t)\right), t \in [T_1, T_2] \quad (4.10)$$

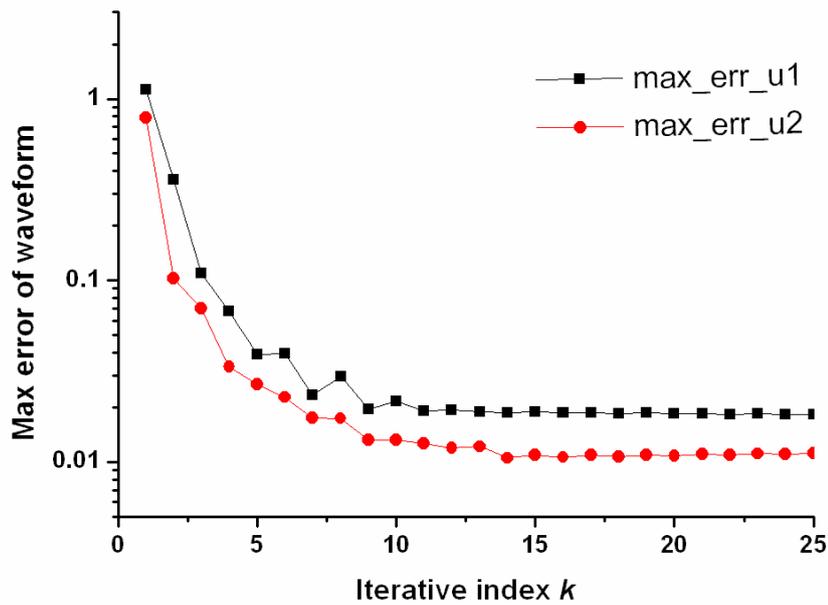

Figure 4. Max error curve of waveforms.

5. Conclusion and Future Work.

WTM is a new distributed algorithm to solve large ODEs. It is based on WR and VTM. If the propagation delays of WTLs are different, WTM would turn to an asynchronous algorithm, similar to DTM [5].

Experiments show that WTM is convergent to solve SPD ODEs. However, its convergence speed and precision are not so impressive, compared to VTM. To speedup WTM, the characteristic impedance waveform $Z(t)$ of WTL should be carefully selected.

Further, WTM might be used to solve the nonlinear ODEs extracted from large scale nonlinear dynamical systems.

Frankly speaking, it is still a doubt whether WTM would be a fast enough algorithm to be accepted by the applied math industry, or just like WR, beautiful but impractical.